\documentclass[11pt]{article}

\usepackage{amscd,amsmath, amssymb, fancyhdr}

\newcommand{\version}{version 4.0,\ \   Sept. 07, 2009}

\numberwithin{equation}{section}

\def\eqref#1{(\ref{#1})}
\newcommand{\goth}{\mathfrak}

\newcommand{\arrow}{{\:\longrightarrow\:}}
\newcommand{\Z}{{\Bbb Z}}
\newcommand{\C}{{\Bbb C}}
\newcommand{\R}{{\Bbb R}}

\def\1{\sqrt{-1}\:}

\newcommand{\calo}{{\cal O}}

\renewcommand{\tilde}{\widetilde}

\renewcommand{\phi}{\varphi}
\renewcommand{\epsilon}{\varepsilon}
\renewcommand{\geq}{\geqslant}
\renewcommand{\leq}{\leqslant}

\newcommand{\End}{\operatorname{End}}

\newcounter{Mycounter}[section]
\newcounter{lemma}[section]
\newcounter{claim}[section]
\renewcommand{\theclaim}{{Claim \thesection.\arabic{claim}}}
\newcommand{\claim}{%
     \setcounter{claim}{\value{Mycounter}}
     \refstepcounter{claim}
     \stepcounter{Mycounter}
     {\noindent \bf \theclaim:\ }}

\newcounter{sublemma}[section]
\newcounter{corollary}[section]
\renewcommand{\thecorollary}{{Corollary \thesection.\arabic{corollary}}}
\newcommand{\corollary}{%
     \setcounter{corollary}{\value{Mycounter}}
     \refstepcounter{corollary}
     \stepcounter{Mycounter}
     {\noindent \bf \thecorollary:\ }}

\newcounter{theorem}[section]
\renewcommand{\thetheorem}{{Theorem \thesection.\arabic{theorem}}}
\newcommand{\theorem}{%
     \setcounter{theorem}{\value{Mycounter}}
     \refstepcounter{theorem}
     \stepcounter{Mycounter}
     {\noindent \bf \thetheorem:\ }}

\newcounter{conjecture}[section]
\newcounter{proposition}[section]
\renewcommand{\theproposition}
       {{Proposition \thesection.\arabic{proposition}}}
\newcommand{\proposition}{%
     \setcounter{proposition}{\value{Mycounter}}
     \refstepcounter{proposition}
     \stepcounter{Mycounter}
     {\noindent \bf \theproposition:\ }}

\newcounter{definition}[section]
\renewcommand{\thedefinition}
       {{Definition~\thesection.\arabic{definition}}}
\newcommand{\definition}{%
     \setcounter{definition}{\value{Mycounter}}
     \refstepcounter{definition}
     \stepcounter{Mycounter}
     {\noindent \bf \thedefinition:\ }}

\newcounter{example}[section]
\newcounter{remark}[section]
\renewcommand{\theremark}{{Remark \thesection.\arabic{remark}}}
\newcommand{\remark}{%
     \setcounter{remark}{\value{Mycounter}}
     \refstepcounter{remark}
     \stepcounter{Mycounter}
     {\noindent \bf \theremark:\ }}

\newcounter{problem}[section]
\newcounter{question}[section]
\makeatletter

\@addtoreset{equation}{section}
\@addtoreset{footnote}{section}
\makeatother

\def\blacksquare{\hbox{\vrule width 5pt height 5pt depth 0pt}}
\def\endproof{\hfill\blacksquare}

\begin{document}
\begin{center}
{\LARGE\bf
Topology of locally conformally K\"ahler\\[2mm]  manifolds with
potential
}\\[3mm]

Liviu Ornea\footnote{Partially supported by a PN II IDEI Grant nr. 529}
and Misha Verbitsky\footnote{Partially supported by the
grant RFBR for support of scientific schools
NSh-3036.2008.2 and RFBR grant 09-01-00242-a.\\


}

\end{center}

{\small
\hspace{0.15\linewidth}
\begin{minipage}[t]{0.7\linewidth}
{\bf Abstract} \\
Locally conformally K\"ahler (LCK) manifolds with
potential are those which admit a K\"ahler covering with a
proper, automorphic, global potential. Existence of a potential
can be characterized cohomologically as vanishing
of a certain cohomology class, called the Bott-Chern class.
Compact LCK manifolds with potential are stable at small
deformations and admit holomorphic embeddings into Hopf
manifolds. This class strictly includes the Vaisman
manifolds. We show that every compact LCK manifold with
potential can be deformed into a Vaisman manifold.
Therefore, every such manifold is diffeomorphic to
a smooth elliptic fibration over a K\"ahler  orbifold.
We show that the pluricanonical condition on LCK manifolds
introduced by G. Kokarev is equivalent to vanishing of
the Bott-Chern class. This gives a simple proof of
some of the results on topology of pluricanonical
LCK-manifolds, discovered by Kokarev and Kotschick.
\end{minipage}
}



\section{Introduction}



The main object of the present paper is
the following notion.

\hfill

\definition
A {\bf  locally conformally K\"ahler  (LCK)} manifold
is a complex Hermitian manifold, with a Hermitian form $\omega$
satisfying $d\omega = \theta\wedge\omega$, where $\theta$ is a closed
 1-form, called {\bf  the Lee form} of $M$.

\hfill

Sometimes an LCK manifold is defined as
a complex manifold which has a K\"ahler covering $\tilde M$, with
the deck transform group acting on $\tilde M$ by conformal homotheties.
This definition is equivalent to the first one, up to
conformal equivalence.

A compact LCK manifold never admits
a K\"ahler structure, unless the cohomology
class $\theta\in H^1(M)$ vanishes (see \cite{Va_80}).
Further on, we shall usually assume that $\theta$
is non-exact, and $\dim_\C M \geq 3$.

LCK manifolds form an interesting class of complex non-K\"ahler
manifolds, including all non-K\"ahler surfaces which are not
class VII. In many situations, the LCK structure becomes useful
for the study of topology and complex geometry of an LCK-manifold.

We would like to investigate the LCK geometry
along the same lines as used to study the
K\"ahler  manifolds. Existence of a K\"ahler structure gives all kinds of
constraints on the topology of $M$ (even-dimensionality of
$H^{odd}(M)$, strong Lefschetz, homotopy formality). It is
thus natural to ask what can we say about the topology of
a compact  LCK manifold.

Not much is known in the general case. We list several known facts:

(1) The LCK manifolds are not necessarily homotopy formal: the
Kodaira surfaces are not homotopy formal (they have
non-vanishing Massey products), but they are LCK (\cite{belgun}).

(2) In \cite{Va_80}, I. Vaisman  conjectured that  for a compact LCK
manifold,
$h^{1}(M)$ should be odd. This was disproven by Oeljeklaus and Toma in
\cite{ot}. In
the same paper Vaisman also conjectured that no compact LCK manifold
can be homotopy equivalent to a K\"ahler manifold. This is still unknown.

(3) All compact Vaisman manifolds (see below) have odd $b_1$.

(4) All non-K\"ahler complex surfaces admit an LCK structure
(\cite{belgun}), except some of Kodaira class VII surfaces.
For Kodaira class VII with $b_2=0$,
LCK structures are known to exist on two types of Inoue
surfaces (Tricerri, \cite{_Tricerri_}), and do not exist on all of the
third type
(Belgun, \cite{belgun}). For $b_2>0$, all Kodaira class VII surfaces
are conjectured to admit a spherical shell
(Ma. Kato, I. Nakamura; see e.g. \cite{_Nakamura:1984_}).
The known examples of minimal Kodaira class VII surfaces
are either hyperbolic or parabolic Inoue surfaces.
LCK structures on  hyperbolic Inoue surfacse were
recently constructed by A. Fujiki and M. Pontecorvo
(\cite{FP}). There are no known examples of non-LCK,
non-K\"ahler surfaces, except the example of
Inoue surface of class $S^+$ considered by
Belgun.

Among the LCK manifolds, a distinguished class is the following:

\hfill

\definition
An LCK manifold $(M,\omega, \theta)$ is called {\bf Vaisman}
if $\nabla\theta=0$, where $\nabla$ is the
Levi-Civita connection of the metric $g(\cdot,\cdot)=\omega(I\cdot,\cdot)$.

\hfill

The universal cover of a Vaisman manifold can be precisely described. We
need the following

\hfill

\definition
A {\bf conical K\"ahler manifold} is a
K\"ahler manifold $(C, \omega)$ equipped with a free, proper
holomorpic flow $\rho:\; \R\times C\arrow C$, with $\rho$ acting
by homotheties as follows: $\rho(t)^*\omega = e^t \omega$. The space of
orbits of $\rho$ is called {\bf a Sasakian manifold}.

\hfill

\theorem\label{_Vaisman_structure_}
(\cite{_ov:Structure_})
A compact Vaisman manifold is conformally equivalent
to a quotient of a conical K\"ahler manifold by $\Z$ freely acting
on $(C, \omega)$ by non-isometric homotheties. Moreover, $M$ admits a
smooth Riemannian
submersion $\sigma:\; M \arrow S^1$, with Sasakian fibers.\hfil \endproof

\hfill

We shall be interested in LCK manifolds whose K\"ahler metric on the
universal cover have global potential. Recall first the following:

\hfill

\definition
Let $(M,I,\omega )$ be a  K\"ahler manifold.
Let $d^c:=-IdI$.
A {\bf K\"ahler potential} is a function
satisfying $dd^c \psi = \omega$. Locally,
a K\"ahler potential always exists, and it
is unique up to adding real parts of holomorphic
functions.

\hfill

\claim
(see e.g. \cite{_Verbitsky_vanishing_})
 Let $(M,\omega, \theta)$
be a Vaisman manifold, and
\[ (\tilde M, \tilde \omega) \stackrel\pi\arrow M\] be
its K\"ahler covering, with $\Gamma\cong \Z$ the
deck transform group: $M = \tilde M/\Gamma$.
Then $\pi^* \theta$ is exact on $\tilde M$:
$\pi^* \theta = d \nu$. Moreover, the function $\psi := e^{-\nu}$
is a K\"ahler potential:  $dd^c \psi =\tilde \omega$. \hfil \endproof

\hfill

\remark
In these notations, let $\gamma\in \Gamma$.
Since $\Gamma$ preserves $\theta$, we have $\gamma^* \nu = \nu + c_\gamma$,
where $c_\gamma$ is a constant. Then $\gamma^* \psi = e^{-c_\gamma}\psi$.
A function which satisfies such a property for
any $\gamma\in \Gamma$ is called {\bf automorphic}.

\hfill

\definition
Let $(M,\omega, \theta)$ be an LCK manifold,
$(\tilde M,\tilde \omega)$ its K\"ahler covering,
$\Gamma$ the deck transform group, $M = \tilde M/\Gamma$,
and $\psi\in C^\infty (\tilde M)$ a K\"ahler potential, $\psi>0$.
Assume that for any $\gamma \in \Gamma$,
$\gamma^*\psi = c_\gamma \psi$, for some constant
$c_\gamma$. Then $\psi$ is called {\bf an automorphic
potential} of $M$.

\hfill

To go on, we need to introduce the following:

\hfill

 \definition
Let  $(M,\omega, \theta)$ be an LCK manifold,
and $L$ a trivial line bundle, associated
to the representation $\mathrm{GL}(2n,\mathbb{R})\ni A \mapsto | \det
A|^{\frac{1}{n}}$, with flat connection
defined as $D:= \nabla_0 + \theta$, where $\nabla_0$
is the trivial connection. Then $L$ is called {\bf the
weight bundle} of $M$. Being flat, its holonomy defines a map $\chi:
\pi_1(M) \arrow \R^{>0}$ whose
image $\Gamma$ is called {\bf the monodromy group of $M$}.

We shall denote with the same letter,
$D$, the corresponding Weyl covariant derivative on $M$.

\hfill

\proposition
(\cite{_ov:MN_})
Let $(M,\omega, \theta)$ be an LCK manifold with an automorphic
potential. Then there exists another LCK metric on $M$
with automorphic potential and monodromy $\Z$. \hfil \endproof

\hfill

\remark
A note on the terminology. In
\cite{_OV:Potential_}, we introduced the
{\bf LCK manifolds with potential}. An LCK manifold is called
{\bf LCK with potential} if it has
an  automorphic potential and monodromy $\Z$.
This terminology can be confusing, because
there exist LCK manifolds with automorphic
potential, without being ``LCK with potential'',
in this sense.

\hfill

The main example of a LCK manifold with potential is the Hopf manifold
$H_A = (\mathbb{C}^n\setminus \{0\})/\langle
A\rangle$, where $A$ is a linear operator with subunitary absolute value
of all eigenvalues, see \cite{_OV:Potential_} where we proved:

\hfill

\theorem
(\cite{_OV:Potential_})
A compact LCK manifold $M$, $\dim_\C M \geq 3$
is LCK with potential
if and only if it admits a holomorphic embedding
to a Hopf manifold $H_A$. $M$ is Vaisman if and only if $A$ is diagonal.
\hfil \endproof

\hfill

The class of compact LCK manifolds is not stable under
small deformations, and the same is true for
Vaisman manifolds (\cite{belgun}).
On the contrary:

\hfill


\theorem
(\cite{_OV:Potential_})
Let $(M,I)$ be a compact complex manifold admitting an LCK
metric with LCK potential. Then any small deformation
of the complex structure $I$ also admits an LCK
metric with LCK potential. \hfil \endproof


\section{Deforming a  compact LCK manifold
with potential to a Vaisman manifold}


\theorem\label{main} Let $(M,\omega,\theta)$,
$\dim_\C M \geq 3$, be an LCK manifold
with potential. Then there exists a small deformation of $M$ which
admits a Vaisman metric.

\hfill

\noindent {\bf Proof:} The idea  is to embed $M$ into a Hopf manifold
defined by a
linear operator $A$, then, using the Jordan-Chevalley decomposition of $A$
to show
that its semisimple part preserves some subvariety of
$\C^n\setminus\{0\}$, thus
yielding an embedding of a small  deformation of $M$ into a  diagonal Hopf
manifold,
the Vaisman metric of which can be pulled back on $M$. We now provide the
details.

\hfill

\noindent{\bf Step 1.} Let $V= \C^n$, $A\in \End(V)$ be an invertible
linear operator
with all eigenvalues $|\alpha_i|<1$, and $H = (V\setminus \{0\})/\langle
A\rangle$ be
the corresponding Hopf manifold, as constructed in \cite{_OV:Potential_}.
One may see that the complex submanifolds
of $H$ are identified with complex subvarieties
$Z$ of $V$, which are smooth outside of $\{0\}$ and are fixed by $A$.
Indeed, by Remmert-Stein theorem (\cite{_Demailly:book_}, chapter II,
\P 8.2), for every complex subvariety $X\subset H$,
the closure of $\pi^{-1}(X)$ is complex analytic in $V=\C^n$,
where $\pi:\; V \backslash 0 \arrow H$ is the natural projection.

\hfill

\noindent{\bf Step 2.}
We are going to prove that any such $Z$ is fixed by the flow
$G_A:=e^{t\log A}$, $t\in \R$, acting on $V$. Let $I_Z$ be
the ideal of $Z$, and let
 $\hat I_Z$ be the corresponding ideal in the completion
of the structural ring $\calo_{V}$ in $\{0\}$. To prove that $I_Z$ is
fixed by $G_A$,
it is enough to show that $\hat I_Z$ is fixed by $G_A$.
However, by definition (see
\cite[Ch. 10]{_Atiyah_MacDonald_}),
$\hat I_Z$ is the inverse
limit
of the projective system $\displaystyle\frac{I_Z}{I_Z \cap {\goth m}^k}$,
where ${\goth m}$ is the maximal ideal of $\{0\}$:
\[ \hat I_Z= \lim\limits_{\leftarrow} \frac{I_Z}{I_Z \cap {\goth m}^k}.
\]
To prove that $\hat I_Z$ is fixed by $G_A$
it only remains to show that $\displaystyle\frac{I_Z}{I_Z \cap {\goth m}^k}$
is fixed by $G_A$. But $\displaystyle\frac{I_Z}{I_Z \cap {\goth m}^k}$
is a subspace in the  vector space
$\calo_{V}/{\goth m}^k$, finite-dimensional by
\cite{_Atiyah_MacDonald_}, Corollary 6.11 and Exercise 8.3,
and such a subspace, if fixed by $A$, is
automatically fixed by $G_A$.

\hfill

\noindent{\bf Step 3.} Now, for any linear operator there
exists a unique decomposition $A:= S U$ in a product
of commuting operators, with $S$ semisimple (diagonal),
and $U$ unipotent, {\em i.e.} its spectrum contains only the number $1$
(this is called
the Jordan-Chevalley decomposition and, in this particular situation, for
operators acting on a finite dimensional vector space
over $\mathbb{C}$, follows easily from the Jordan canonical form).
Consequently, for any finite-dimensional representation of $GL(n)$,
any vector subspace which is fixed by $A$, is also fixed by $S$.
By the argument in Step 2, this proves that  $S$
fixes the ideal $\hat I_Z$, and the subvariety
$Z\subset V$.

\hfill

\noindent{\bf Step 4.} The diagonal Hopf variety
$H_S:= (V\setminus \{ 0\})/\langle S\rangle$ contains a Vaisman
submanifold $M_1 := (Z\setminus \{0\})/\langle S\rangle$. Since $S$ is
contained in
a closure of a $GL(V)$-orbit of $A$, we have also shown that
$M_1$ can be obtained as an arbitrary small deformation
of $M$.
\endproof

\hfill

 It is well known that every compact Vaisman manifold
 is diffeomorphic to a quasiregular Vaisman manifold, \cite{_ov:Immersion_}.
By definition, the latter is an elliptic fibration over a
K\"ahler  orbifold. As a
 consequence, the above result gives us  the possibility
 to obtain topological information about LCK manifolds
 with potential from the topology of projective
 orbifolds. But, once a compact LCK manifold admits an
 automorphic potential on a covering, it can be deformed
 to a LCK manifold with a proper potential (\cite{_ov:MN_},
Corollary 5.3). Hence we have:

\hfill

\corollary \label{fund_group}
The fundamental group of a
compact LCK manifold $M$ with an automorphic potential
admits an exact sequence
\[
0 \arrow G \arrow \pi_1(M) \arrow \pi_1(X) \arrow 0
\]
where $\pi_1(X)$ is the fundamental group of a K\"ahler  orbifold,
and $G$ is a quotient of $\Z^2$ by a subgroup of rank $\leq 1$.

\hfill

\noindent{\bf Proof:} Replacing $M$ with a diffeomorphic
Vaisman manifold, we may assume that $M$ is
a quasiregular Vaisman manifold, elliptically fibered
over a base $X$. The long exact sequence of homotopy gives
\[
\pi_2(X) \stackrel \delta
\arrow \pi_1(T^2) \arrow \pi_1(M) \arrow \pi_1(X) \arrow 0
\]
The boundary operator $\delta$ can be described
as  follows. Let \[ \gamma:\; \Z^2 \arrow H^2(X)\]
be the map  representing the Chern
classes of the corresponding $S^1\times S^1$-fibration.
We may interpret this map as a differential
of the corresponding Leray spectral sequence, which
gives us an exact sequence
\[
0 \arrow H^1(X) \arrow H^1(M) \arrow
H^1(T^2)\stackrel \gamma\arrow H^2(X).
\]
Dualizing and using the Hurewicz theorem,
we obtain that the  boundary map $\pi_2(X) \stackrel \delta
\arrow \pi_1(T^2)$ is obtained as a
composition of $\gamma^*$ and the
Hurewicz homomorphism $\pi_2(X) \arrow H^2(X)$.
The Chern classes of the $S^1\times S^1$-fibration
 are easy to compute: one of them
is trivial (because $M$ is fibered over a circle),
and the other one is non-trivial,
because $M$ is non-K\"ahler, and the total space of
an isotrivial elliptic fibration with trivial Chern classes is K\"ahler.
Therefore, the image of $\delta$ has rank $\leq 1$ in $\pi_1(T^2)$.
\endproof


\section{Pluricanonical LCK manifolds are diffeomorphic
to Vaisman manifolds}


In \cite{_Kokarev:plcan_},
G. Kokarev introduced the following notion:

\hfill

\definition (\cite{_Kokarev:plcan_})
Let $(M, \omega,\theta)$ be an LCK manifold. Then
$M$ is called {\bf pluricanonical} if
$(\nabla\theta)^{1,1}=0$, where $(\cdot)^{1,1}$
denotes the $I$-invariant part of the tensor.

\hfill

For this class of LCK manifolds, Kokarev and Kotschick
generalized an important result by Siu and Beauville
\cite{_Siu:rigidity_} from K\"ahler geometry:

\hfill

\theorem\label{surj_hom}
(\cite{koka_kot})
Let $M$ be a compact pluricanonical LCK manifold,
such that $\pi_1(M)$ admits a surjective
homomorphism to a non-abelian free group.
Then $M$ admits a surjective holomorphic
map with connected fibers to a compact
Riemannian surface.

\hfill

Clearly, the pluricanonical condition is weaker than the
Vaisman one. But in the cited papers, no other examples of pluricanonical
LCK manifolds
are provided but Vaisman ones.

We now prove that the pluricanonical condition is
equivalent with the existence of an automorphic potential
on a K\"ahler covering.

Indeed, the Levi-Civita connection
$\nabla$ and the Weyl connection $D$ on $M$ are related by
the formula (\cite{drag}):
$$
\nabla -D = \frac 12
(\theta\otimes\mathrm{id}+\mathrm{id}\otimes\theta-g\otimes
\theta^\sharp).
$$
Applied on $\theta$, this gives:
$$ \nabla\theta-D\theta=-\theta\otimes\theta +\frac 12 g.$$
Hence, the pluricanonical condition $(\nabla\theta)^{1,1}=0$ is translated
into
\[
(D\theta)^{1,1} = (\theta\otimes \theta)^{1,1}-\frac 12 g.
\]
Since $D$ is torsion-free,
this is equivalent to
$$d(I\theta) = \omega - \theta\wedge I\theta.$$
 But we
can prove
the following:

\hfill

\claim
Let $(M, \omega, \theta)$ be an LCK manifold, and
$\theta^c:= I(\theta)$ the complex conjugate of the Lee form.
Then the condition
$d\theta^c = \omega - \theta\wedge \theta^c$
is equivalent to the existence of an automorphic potential on a K\"ahler
covering on
which the pull-back of $\theta$ is exact.

\hfill

\noindent{\bf Proof:}
Let $\tilde M$ be a covering of $M$ on which the pull-back
of $\theta$ is exact. Denote, for convenience, with the
same letters the pull-backs to $\tilde M$ of $\theta$,
$\omega$ and $D$. Observe that the Levi-Civita
connection of the
K\"ahler metric on $\tilde M$ globally conformal with
$\omega$ is precisely $D$.
Let $\psi:= e^{-\nu}$, where $d\nu = \theta$.
Then
\[  d d^c  \psi=
   -e^{-\nu} d d^c \nu + e^{-\nu} d\nu \wedge d^c \nu
   =  e^{-\nu} (d^c\theta + \theta\wedge I\theta)= \psi \omega,
\]
and hence the pluricanonical condition implies that $\psi$
is an automorphic potential for the K\"ahler metric $\psi\omega$. The
converse is true
by the same argument.
\endproof

As we know that the existence of an automorphic potential
on a covering allows the deformation to a LCK manifold
with potential, taking into account \ref{main}, we obtain
the following

\hfill

\corollary Any compact pluricanonical LCK manifold is
diffeomorphic to a Vaisman manifold.

\hfill

In view of \ref{fund_group}, the restrictions on
the fundamental group of a pluricanonical LCK manifold
obtained in \cite{koka_kot} using a generalization of
Siu's arguments for harmonic maps can be directly
obtained for LCK manifolds which admit an automorphic
potential on a K\"ahler covering, by using \ref{fund_group} and
the corresponding results for K\"ahler manifolds.
For instance, the following result can be proven
(this is one of two cases Corollary 3.3 of \cite{koka_kot}).

\hfill

\claim
A non-abelian free group cannot be the
fundamental group of an LCK manifold $M$
admitting an LCK potential.

\hfill

\noindent{\bf Proof:} As  follows from \ref{fund_group},
$\pi_1(M)$ fits into an exact sequence
\[
0 \arrow G \arrow \pi_1(M) \arrow \pi_1(X) \arrow 0,
\]
where $G$ is an abelian group of rank $\geq 1$. Since all
non-trivial normal subgroups of infinite index
in a free group are free, by Nielsen-Schreier theorem, and
infinitely generated (see e.g. \cite{_Greenberg:free_}),
$G$ cannot be a normal subgroup of a free group,
unless $G$ is trivial.
\endproof

\hfill

\hfill

\noindent{\bf Acknowledgements:} We are grateful to Victor
Vuletescu for many useful remarks and for pointing an error
in an earlier version of this paper. We are also grateful to the anonymous
referees for their comments
which improved the presentation.

{\small

\noindent {\sc Liviu Ornea\\
University of Bucharest, Faculty of Mathematics, \\14
Academiei str., 70109 Bucharest, Romania. \emph{and}\\
Institute of Mathematics ``Simion Stoilow" of the Romanian Academy,\\
21, Calea Grivitei Street
010702-Bucharest, Romania }\\
\tt Liviu.Ornea@imar.ro, \ \ lornea@gta.math.unibuc.ro

\hfill

\noindent {\sc Misha Verbitsky\\
{\sc  Institute of Theoretical and
Experimental Physics \\
B. Cheremushkinskaya, 25, Moscow, 117259, Russia }\\
\tt verbit@maths.gla.ac.uk, \ \  verbit@mccme.ru
}
}

\end{document}